\documentclass[a4paper,11pt]{article}
\usepackage[plainpages=false]{hyperref}
\usepackage{amsfonts,latexsym,rawfonts,amsmath,amssymb,amsthm,mathrsfs}
\usepackage{amsmath,amssymb,amsfonts,latexsym,lscape,rawfonts}

\usepackage[all]{xy}
\usepackage{eufrak}
\usepackage{makeidx}         % allows index generation
\usepackage{graphicx,psfrag}

\usepackage{array,tabularx}

\usepackage{setspace}

\newtheorem{thm}{Theorem}[section]
\newtheorem{cor}[thm]{Corollary}
\newtheorem{lem}[thm]{Lemma}

\newtheorem{prop}[thm]{Proposition}
\newtheorem{q}[thm]{Problem}

\theoremstyle{remark}
\newtheorem{rmk}[thm]{Remark}

\theoremstyle{definition}

\def \N {\mathbb N}
\def \C {\mathbb C}

\def \R {\mathbb R}

\def \P {\mathbb P}
\def \O {\mathcal O}

\def \B {\mathcal B}

\def \H {\mathcal H}
\def \L {\mathcal L}
\def \p {\partial}

\def \fid {\dot{\phi}}
\def \fidd {\ddot{\phi}}
\def \fiddd {\dddot{\phi}}
\def \pid {\dot{\phi_k}}
\def \pidd {\ddot{\phi_k}}
\def \piddd {\dddot{\phi_k}}
\def \hd {\dot{H}}
\def \hdd {\ddot{H}}
\def \hddd {\dddot{H}}

\begin{document}
\title{\textbf{Space of K\"ahler metrics (V)---\\
K\"ahler quantization}}
\author{Xiuxiong Chen\footnote{Partially supported by  NSF grant.} and Song Sun}
\date{}

\maketitle

\begin{abstract}
Given a polarized K\"ahler manifold $(X, L)$. The space $\H$ of
K\"ahler metrics in $2\pi c_1(L)$ is an infinite dimensional
Riemannian symmetric space. As a metric space, it has non-positive
curvature. There is associated to $\H$ a sequence of finite
dimensional symmetric spaces $\B_k$($k\in \N$) of non-compact
type. We prove that $\H$ is the limit of $\B_k$ as metric spaces
in the weak sense.  As applications, this provides more geometric
proofs of certain known geometric properties of the space $\H$.

\end{abstract}

\section{Introduction}
Let $(X, \omega, J)$ be an $n $ dimensional K\"ahler manifold. By
\cite{D1} this gives rise to two infinite dimensional symmetric
spaces: the Hamiltonian diffeomorphism group $Ham(M, \omega)$ and
the space $\H$ of smooth K\"ahler potentials  under the natural
(Weil-Petersson type) $L^2$ metric. The relation between these two
spaces is analogous to that between a finite dimensional compact
group $G$ and its noncompact dual $G^{\C}/G$, where $G^{\C}$ is a
complexification of $G$. As in the finite dimensional case, at
least formally, $Ham(M, \omega)$ has non-negative sectional
curvature, while $\H$ has non-positive sectional curvature. It is
proved by E. Calabi and the first author in \cite{CC} that $\H$ is
non-positively curved in the sense of Alexandrov. On the other
hand, it is well-known in the literature of geometric quantization
that  $\H$ is the limit of a sequence of finite dimensional
symmetric spaces $\B_k(k\in\N)$ when $X$ is (K\"ahler) polarized.
Indeed, the polarization is given by an ample line bundle $L$ over
$X$, and we consider the space $H^0(X, L^k)$ of holomorphic
sections of $L^k$ for large enough $k$. Here $k^{-1}$ plays the
role of Planck constant $\hbar$, while $k\rightarrow\infty$ should
correspond to the process of taking the classical limit. Denote by
$N_k$ the dimension of $H^0(X, L^k)$. By the Riemann-Roch theorem
this is a polynomial in $k$ of degree $n$ for sufficiently large
$k$. Then $\B_k$ is chosen to be the space of positive definite
Hermitian forms on $H^0(X, L^k)$. It could also be viewed as the
symmetric space $GL(N_k;\C)/U(N_k)$. The spaces $\B_k$ are related
to $\H $ through two naturally defined maps:
$Hilb_k:\H\rightarrow\B_k$ and $FS_k:\B_k\rightarrow\H$.(The
precise definition of these maps will be given in the next
section). In \cite{T}, G. Tian proved that given any $\phi\in \H$,
$FS_k\circ Hilb_k(\phi)\rightarrow\phi$ in $C^4$ topology. Using
Tian's  peak section method and canonical coordinates, W-D Ruan
proved the $C^{\infty}$ convergence in \cite{Ruan}. In \cite{Z} S.
Zelditch beautifully generalized Tian's theorem and derived
$C^{\infty}$ convergence from the asymptotic expansion of the
Bergman kernel. If we denote by $\H_k$ the image of $FS_k\circ
Hilb_k$, then these results tell us that
$$\H=\overline{\bigcup_{k\in\N}\H_k},$$
where the closure is taken in $C^{\infty}$ topology.\\

 In \cite{D2}, S. K. Donaldson suggests that the geometry of $\B_k$
should also converge to that of $\H$ because of the previous
strong similarity between the two. In particular, the geodesics in
$\H$ should be approximated by geodesics in $\H_k$.  In \cite{C1},
the first author  proved the existence of $C^{1,1}$ geodesics and
it subsequently lead to  many interesting applications in K\"ahler
geometry (see \cite{CT2}, \cite{C2} for further references).
 The limitation of $C^{1,1}$ regularities is purely technical.
In a very interesting paper(\cite{PS1}), Phong-Sturm proved that
the $C^{1,1}$ geodesics in $\cal H$ are the weak $C^0$ limits of
Bergman geodesics, assuming the existence of $C^{1,1}$ geodesics.
It would provide a canonical smooth approximation of $C^{1,1}$
geodesics if one could show the convergence is in $C^{1,1}$
topology. More evidence comes from a series of beautiful works by
S.Zeltdich and his collaborators, see Song-Zelditch(\cite{SZ}),
Rubinstein-Zelditch(\cite{RZ}). They proved that on toric
varieties both geodesics and harmonic maps in $\H$(automatically
smooth) are $C^2$ limits of corresponding objects in $\H_k$.
Recently, J. Fine (\cite{Fi}) proved the remarkable result that
the Calabi flow in $\H$ could be approximated by balancing flows
in $\B_k$.  In this paper, we shall prove the following
convergence of
geodesic distance:\\

\begin{thm}
\label{main} Given any $\phi_0, \phi_1\in H$, we have
$$\lim_{k\rightarrow\infty}k^{-\frac{n+2}{2}}d_{\B_k}(Hilb_k(\phi_0), Hilb_k(\phi_1))= d_{
\H}(\phi_0, \phi_1).\\$$
\end{thm}

\begin{rmk}
It is easy to identify the scaling factor by simply taking
$\phi_1=\phi_0+c$, then $d(\phi_0,\phi_1)=c$, while
$d(Hilb_k(\phi_0), Hilb_k(\phi_1))=c\cdot k\sqrt{N_k}$. This
scaling also indicates that the limit can not have bounded
curvature, which can also be speculated from the expression of the
infinitesimal curvature of $\H$:
$$R_{\phi}(\phi_1, \phi_2)=-\frac{1}{4}\frac{||\{\phi_1, \phi_2\}||_{L^2}^2}{||\phi_1||_{L^2}^2||\phi_2||_{L^2}^2}.$$
\end{rmk}

\begin{rmk}
This theorem indicates that the $L^2$ metric defined in \cite{D1}
is in a sense canonical.  We notice that there are also natural
affine structures on both $\H$ and $\B_k$, given by embedding them
as convex subsets of the affine spaces $C^{\infty}(X)$ and the
space of all Hermitian forms on $H^0(X, L^k)$ respectively. It
would be interesting to see whether the affine structures on $\H$
and $\B_k$ have nice relations with each other
as the above hyperbolic geometric structures do. \\
\end{rmk}

 From the proof it follows that the convergence is uniform if
both potentials vary in a $C^l$ compact neighborhood for large
$l$. So an easy corollary is the
convergence of angles:\\

\begin{cor}
Given three points $\phi_i(i=1,2,3)$ in $\H$, let
$H_{k,i}=Hilb_k(\phi_1)$.  Then
$$\lim_{k\rightarrow\infty}\angle H_{k,1}H_{k,2}H_{k,3}=\angle \phi_1\phi_2\phi_3.$$
\end{cor}

This corollary leads to that $\H$ is non-positively curved in the
sense of Alexandrov, as was originally proved in \cite{CC}.\\

 Theorem \ref{main} together with theorem 5 in \cite{Fi}
imply the following corollary, because the downward gradient flow
of a geodesically convex function on a finite dimensional manifold
is distance decreasing:

\begin{cor}\label{Calabi flow} (Calabi-Chen\ \textbf{\cite{CC}} \textbf{\cite{C2}}) Calabi flow in $\H$ decreases geodesic distance as
long as it exits.
\end{cor}

 The ``quantization" approach, namely using the approximation by $\B_k$ to
handle problems about the K\"ahler
 geometry of $(X, \omega, J)$, turns out to be quite powerful and intriguing.
It was shown in \cite{D1} that the existence and
 uniqueness problems in K\"ahler geometry are related to the
 geometry of $\H$. More precisely, the uniqueness of extremal
 metrics is implied by the existence of smooth geodesics
 connecting any two points in $\H$, and the non-existence is
 conjectured to be equivalent to the existence of a geodesic ray
 in $\H$ where the K-energy is strictly decreasing near infinity.
 The technical problem comes from the lack of regularity of
 geodesics in $\H$. There are two ways
at present to circumvent this problem. The first way is to try to
study the geodesic equation directly. One can use the continuity
method, as in \cite{C1}. In \cite{C1},the first author constructed
a continuous family of $\epsilon$-approximate geodesics converging
in weak $C^{1,1}$ topology to a $C^{1,1}$ geodesic. This gives the
proof of the uniqueness of extremal metrics when $c_{1}(X) \leq0$
and the other conjecture of Donaldson that $\cal H$ is a metric
space. In \cite{CT2} a new partial regularity was derived from
studying the complex Monge-Amp\`ere equation associated to the
geodesic equation. Solutions with such regularity already have
fruitful applications. The second way is to exploit the
approximation of $\H$ by $\B_k$ This requires a polarization. The
essential thing is to try to quantize the whole problem.  Many
problems have been settled using either of the two approaches: the
uniqueness of cscK metrics (\cite{CT2}, \cite{D1}), the lower
bound of Mabuchi energy assuming the existence of cscK metrics(
\cite{CT2}, \cite{D5}), the lower bound of Calabi
energy (\cite{C3}, \cite{D6}), etc.\\

Due to \cite{D5}, the Mabuchi functional (K-energy) $E$ has a nice
quantization defined on $\B_k$, which we denote by $Z_k$. More
precisely, up to a constant, given any $\phi\in\H$, we have
$Z_k(Hilb_k(\phi))\rightarrow E(\phi)$. It is well known that $E$
is convex along smooth geodesics in $\H$, but not known to be
convex along $C^{1,1}$ geodesics, although it is well-defined
(c.f. \cite{C4}). One good thing is that $Z_k$ is convex on a
finite dimensional space $\B_k$, which is geodesically complete.
This makes it possible to derive some weak convexity of $E$, such
as an alternative proof of two inequalities of the first author. Namely\\

\begin{cor}\label{inequality 2} (\textbf{\cite{C2}}) For any $\phi_0,
\phi_1\in\H$, let $\phi(t)(t\in[0,1])$ be the $C^{1,1}$ geodesic
connecting them. Then $$dE_{\phi_0}(\fid(0))\leq
dE_{\phi_1}(\fid(1)).$$
\end{cor}

\begin{cor} \label{inequality 1} (\textbf{\cite{C2}}) For any $\phi_0,
\phi_1\in\H$, we have
$$E(\phi_1)-E(\phi_0)\leq d(\phi_0, \phi_1)\cdot
\sqrt{Ca(\phi_1)}.$$
\end{cor}

Both corollaries are important in the study of K\"ahler geometry
through $\H$.  By Theorem 3.1 in \cite{CC}, corollary
\ref{inequality 2} immediately gives corollary \ref{Calabi flow}.
Using these corollaries, the first author proved the sharp lower
bound of the Calabi energy in any given K\"ahler class. The
original proof of these two corollaries depends heavily on the
delicate regularity results of Chen-Tian \cite{CT2}. So at least
in the algebraic case, one would like to see
 a more geometric proof of these.\\

\noindent {\bf Organization:} In section 2 we briefly recall the
geometry of the space $\H$ of K\"ahler potentials and the related
finite dimensional spaces $\B_k$. In section 3 we prove the
convergence of geodesic distance, using the existence of $C^{1,1}$
geodesics in $\H$, the non-positivity of curvature of $\B_k$, and
the convergence of infinitesimal geometry. In section 4 the weak
convexity of K-energy is discussed and we prove corollary
\ref{inequality 2} and corollary \ref{inequality 1}. The idea is
that the K-energy could be approximated by convex functions $Z_k$
on $\B_k$, while the gradient of $Z_k$ approximates the Calabi
functional, which is the norm of  gradient of K-energy. So
corollary \ref{inequality 1} follows. To prove corollary
\ref{inequality 2}, we need to estimate the difference between the
initial direction of an ``almost" geodesic $\gamma$(i.e
$|\ddot{\gamma}|$ is small) in $\B_k$ and that of the genuine
geodesic connecting the two end points $\gamma(0)$ and
$\gamma(1)$. This is done in lemma \ref{Riemannian lemma 2}. In the last section, we discuss some open problems.\\

\noindent {\bf Acknowledgments:} Both named authors wish to thank
Professor S. Zelditch for his interest in this work. The second
author wishes to thank Professor S. Donaldson, Professor B. Lawson
and Professor G. Tian for their helpful comments. He also would
like to thank Garrett Alston for
suggestions improving the exposition of this article.\\

\section{Preliminaries} Let $(X, \omega, J)$ be a K\"ahler manifold.
Here we always assume it is polarized, i.e. there is a holomorphic
line bundle $L$ whose first Chern class is $[\omega]/2\pi$. Fix
the base metric $\omega$ and define the space of K\"ahler
potentials as
$$\H=\{\phi\in
C^{\infty}(M)|\omega+\sqrt{-1}\partial\bar{\partial}\phi>0\}.$$
$\H$ can also be identified with the space of Hermitian metrics on
$L$ with positive curvature form. In the following we will not
distinguish between these two meanings. We write $\omega_{h}$ to
denote $\sqrt{-1}$ times the curvature of $h$ for a Hermitian
metric $h$, and $\omega=\omega_{h_0}$. Thus,
$$\omega_{\phi}:=\omega_{e^{-\phi}h_0}=\omega+\sqrt{-1}\partial\bar{\partial}\phi.$$
Denote
$$d\mu_{h}=d\mu_{\phi}=\frac{\omega_{h}^n}{(2\pi)^n n!}.$$
There is a Weil-Petersson type metric defined on $\H$ by
$$(\delta_1\phi, \delta_2\phi)_{\phi}=\int_X \delta_1\phi\delta_2\phi d\mu_{\phi},$$
for any $\delta_1\phi,\delta_2\phi\in T_{\phi}\H=C^{\infty}(X)$.
Due to \cite{D1}, $\H$ is formally an infinite dimensional
symmetric space, and the sectional curvature of the Levi-Civita
connection is given by:
$$R(\delta_1\phi, \delta_2\phi)=-\frac{1}{4}||\{\delta_1\phi,
\delta_2\phi\}_{\phi}||^2.$$ It was confirmed in \cite{CC} that
$\H$ satisfies the triangle comparison theorem of a non-positive
Alexandrov space, by using so-called $\epsilon$-approximate
geodesics. The geodesic equation in $\H$ is
\begin{equation}
\fidd-|\nabla_{\phi}\fid|_{\phi}^2=0, \label{G1} \end{equation}
where we have used the notation of complex gradient.
The following theorem is proved in \cite{C1} by the first author:\\

\begin{lem} \label{geodesic approximation lemma}
(Geodesic Approximation Lemma) Given any $\phi_0, \phi_1\in\H$,
there is a positive number $\epsilon_0$, and a one parameter
smooth family of smooth curves $\phi_{\epsilon}(\cdot):
[0,1]\rightarrow\H \ (\epsilon\in (0, \epsilon_0])$, such that the
following holds:\\

(1). For any $\epsilon\in (0, \epsilon_0]$, $\phi_{\epsilon}$ is
an $\epsilon$-approximate geodesic, i.e. it solves the following
equation:
$$(\fidd_{\epsilon}-|\nabla\fid_{\epsilon}|^2)\omega_{\phi_{\epsilon}}^n=\epsilon\cdot\omega^n, $$
and $\phi_{\epsilon}(0)=\phi_0$, $\phi_{\epsilon}(1)=\phi_1$.\\

(2). There exists a $C>0$, such that for all $\epsilon\in (0,1]$
and $t\in [0,1]$, we have
$$|\phi_{\epsilon}(t)|+|\fid_{\epsilon}(t)|\leq C,$$ and
$$|\fidd_{\epsilon}(t)|\leq C.$$

(3). $\phi_{\epsilon}(\cdot)$ converges to the unique $C^{1,1}$
geodesic connecting $\phi_0$ and $\phi_1$
 in the weak $C^{1,1}$ topology when $\epsilon
 \rightarrow 0$. \\

(4). We have uniform estimates when $\phi_0$ and $\phi_1$
varies in a $C^k$ compact set for some large $k$. \\
\end{lem}
The space $\H$ has a global flat direction given by addition of a
constant. The de Rham decomposition theorem turns out to be true
in this case. There is a functional $I$ giving the isometric
decomposition
$$\H=\H_0 \oplus \R.$$ $I$ is only defined up to an addition of
a constant, and its derivative is given by:
$$\delta I= \int_X \delta\phi d\mu_{\phi}.$$
 So $I$ is linear
along $C^{1,1}$ geodesics in $\H$. $\H_0$ is then an arbitrary
level set of $I$, and it can also be regarded as the space of
K\"ahler metrics cohomologous to $\omega$. An interesting thing to
notice is that $\H$ is embedded as a convex subset into an affine
linear space $C^{\infty}(M)$. It is easy to show that $I$
is indeed convex along any linear path in $\H$.\\

 There is a well known K-energy functional
$E$ on $\H$, defined up to an additive constant, with its
variation given by: \begin{equation}\delta E=-\int_X
(S-\underline{S})\delta{\phi}d\mu_{\phi}.
\label{Kenergy}\end{equation} Here $\delta{\phi}$ is the
infinitesimal variation of $\phi$, and $S$ is the scalar curvature
of $\omega_{\phi}$. From the point of view of \cite{D0}, this is a
natural convex functional associated to a moment map of a
Hamiltonian action (for infinite dimensional manifold). Indeed, by
a straightforward calculation, if $\phi(t)$ is a smooth geodesic
in $\H$, then
$$\frac{d^2}{dt^2}E(\phi(t))=\int_X |\mathcal D\fid|^2d\mu_{\phi}\geq 0,$$
where $\mathcal D$ is the Lichnerowicz Laplacian operator. By
studying the explicit expression of $E$, it is shown in \cite{C4}
that  $E$ is well-defined for $C^{1,1}$ K\"ahler potentials, but
it is not obvious from the definition
that $E$ is convex along $C^{1,1}$ geodesics.\\

There is a sequence of finite dimensional symmetric spaces $\B_k$
consisting of all positive definite Hermitian forms on $H^0(X,
L^k)$. This can be identified with $GL(N_k;\C)/U(N_k)$ by choosing
a base point in $\B_k$, where $N_k=dim H^0(X, L^k)$. We shall
review some basic facts about these symmetric spaces. $U(N)$ is a
compact Lie group and admits a natural bi-invariant Riemannian
metric given by $$(A_1, A_2)_A= -Tr(A_1A^{-1}\cdot A_2A^{-1}),
$$for any $A_1, A_2\in T_A U(N)$. The sectional curvature is given
by:
$$R(A_1, A_2)=\frac{1}{4}||[A_1, A_2]||^2.$$ The geodesic equation
is $$\ddot{A}=\dot{A}A^{-1}\dot{A}.$$ It is then clear that the
geodesics all come from one parameter subgroups given by the usual
exponential map. Now the non-compact dual of $U(N)$ is
$GL(N;\C)/U(N)$--the space of positive definite $N\times N$
Hermitian matrices. We can explicitly  write down the metric on
it. For any $H_1, H_2\in T_H (GL(N;\C)/U(N))$, define
$$(H_1, H_2)_H=Tr(H_1H^{-1}\cdot H_2 H^{-1}).$$
Then the sectional curvature becomes non-positive:
$$R(H_1, H_2)=-\frac{1}{4}||[H_1, H_2]_H||^2.$$
A path $H(t)$ is a geodesic if it satisfies the following
equation:
\begin{equation}\ddot{H}=\dot{H}H^{-1}\dot{H}. \label{G2}\end{equation} This
looks very similar to equation $(\ref{G1})$, except that it is a
nonlinear ODE instead of a nonlinear PDE. In our setting, we get a
sequence of equations $(\ref{G2})$ depending on $k$. For this
reason, we may say that this sequence of equations ``quantizing"
equation $(\ref{G1})$. It is easy to see that all the geodesics in
this space are of the form
$$H(t)=H(0)^{\frac{1}{2}}\exp(tA)H(0)^{\frac{1}{2}},$$ for some initial point $H(0)$ and
initial tangent vector $A$.\\

Similar to the infinite dimensional case, $\B_k$ also admits an
isometric splitting given by the function $I_k=Log det$, defined
up to an additive constant. It is easy to see that $I_k$ is linear
along geodesics in $\B_k$ and convex along a linear path when we
regard $\B_k$ as a convex subset of the affine linear space of all
$N_k\times N_k$ complex matrices. The splitting corresponds to the
following isometry:
$$GL(N;\C)/U(N)\simeq SL(N;\C)/SU(N)\times \R.$$
We will see in   section 4 that $I_k$ indeed ``quantizes" the
$I$ functional.\\

To ``quantize " functionals on $\H$, there are two natural maps
relating $\H$ and $\B_k$, where we have adopted the notation of
\cite{D5}:
$$Hilb_k:\H\rightarrow \B_k;$$
$$FS_k: \B_k\rightarrow \H.$$
Explicitly, given $h\in \H$, and $s\in H^0(X, L^k)$, we define
$$||s||_{Hilb_k(h)}^2=\int_X |s|^2_{h}d\mu_h.$$ For $H\in \B_k$,
pick an orthonormal basis $\{s_{\alpha}\}$ of $H^0(X, L^k)$ with
respect to $H$. Define
$$FS_k(H)=\frac{1}{k}log\sum_{\alpha}|s_{\alpha}|^2_{h^k_0},$$
where $h_0$ is the base metric in $\H$. In particular, we have
$$\sum_{\alpha}|s_\alpha|^2_{e^{-k\cdot FS_k(H)}h_0}\equiv1.$$
Notice that $\omega_{FS_k(H)}$ coincides with the pull back of the
Fubini-Study metric on $\C \P^{N_k-1}$ under the projective
embedding induced by $\{s_{\alpha}\}$, while $FS_k(H)$ coincides
with the induced metric on the pull back of $\O(1)$. Let $\H_k$ be
the image of $FS_k$, then the well-known Tian-Yau-Zelditch theorem
says that
$$\overline{\bigcup_k\H_k}=\H.$$ To be more precise, for any
$\phi\in\H$, and $k>0$, we choose an orthonormal basis
$\{s_{\alpha}\}$ of $H^0(X, L^k)$ with respect to $Hilb_k(\phi)$,
and define the density of state function:
$$\rho_k(\phi)=\sum_{\alpha}|s_{\alpha}|^2_{h^k}.$$
Then we have the following  $C^{\infty}$ expansion:\\

\begin{lem} \label{tianyauzelditch}(\cite{Z}, \cite{Lu1})
\begin{equation}\rho_k(\phi)=k^n+
A_1(\phi)k^{n-1}+A_2(\phi)k^{n-2}+\cdots,
\label{TYZ}\end{equation} with $A_1(\phi)=\frac{1}{2}S(\phi)$,
where $S$ denotes the scalar curvature. Moreover, the expansion is
uniform in that for any $l$ and $R\in \N$,
$$||\rho_k(\phi)-\sum_{j\leq R}A_j k^{n-j}||_{C^l}\leq C_{R, l}k^{n-R},$$ where $C_{R,l}$ only depends on $R$ and $l$.\\
\end{lem}

\begin{rmk} \label{Lu} Later we will need a generalization of the
previous expansion theorem. We want to differentiate the expansion
of the density of state function several times along a smooth path
$\phi(t)$ in $\H$, and still to get a uniform expansion.
 This could be done following the arguments of Z. Lu, see \cite{LT}. \\
\end{rmk}

 Now if we let $\phi_k= FS_k\circ Hilb_k
(\phi)$, then
$$\phi_k-\phi=
\frac{1}{k}Log \rho_k(\phi)\rightarrow 0, $$ as
$k\rightarrow\infty$. \\

To quantize the K-energy, following \cite{D5}, we define
$$\L_k=I_k\circ Hilb_k+\frac{kd_k}{V}I,$$ and
$$Z_k=\frac{kd_k}{V}I\circ FS_k+I_k+d_k(Log V-Log d_k),$$
where $d_k=dim H^0(X, L^k)$, and $V=Vol_{[\omega]}(X)$. Notice
that the definition of both functionals requires choices of  base
points in both $\H$ and $\B_k$. From now on, we always fix the
same base points for defining these two functionals. By a
straightforward calculation,
$$\delta\L_k=\int[\Delta\rho_k-k\rho_k]_{\phi}\delta\phi
d\mu_{\phi},$$ where we define for any function $f$
$[f]_{\phi}:=f-\frac{1}{V}\int f d\mu_{\phi}=f-\underline{f}$.
From (\ref{TYZ}), we see that
$$[\Delta\rho_k-k\rho_k]_{\phi}\rightarrow -\frac{1}{2}k^n(S-\underline{S}). $$
Thus, there are constants $c_k$, such that
$$\frac{2}{k^n}\L_k+c_k\rightarrow E,$$
where the convergence is uniform on $C^{l}(l\gg1)$ bounded subsets
in $\H$. From now on, we will denote the left hand side by $\L_k$,
and $Z_k+c_k$ by $Z_k$. \\

In \cite{D5} the following relation between $\L_k$ and $Z_k$ was
shown:
\begin{eqnarray*}
\L_k(FS_k(H))-Z_k(H)&=& Log det (Hilb_k \circ FS_k(H))-Log det
H-d_k Log\frac{V}{d_k}\leq 0;
\end{eqnarray*}
and \begin{eqnarray*} Z_k(Hilb_k(\phi))-\L_k(\phi)&=&
\frac{kd_k}{V}(I(FS_k\circ Hilb_k(\phi))-I(\phi))+d_k
Log\frac{V}{d_k}\leq 0.
 \end{eqnarray*}
Thus, given any $\phi\in \H$, we have
$$\L_k(\phi_k)=\L_k(FS_k \circ Hilb_k(\phi))\leq Z_k(Hilb_k(\phi))\leq \L_k(\phi).$$
Since $\phi_k$ converges to $\phi$ smoothly, and $\L_k$ converges
to $E$ uniformly on $C^l(l>>1)$ bounded subsets in $\H$, we
immediately
obtain:\\

\begin{lem} \label{ZquantizesE}$Z_k$ quantizes $E$  in the sense that given any
$\phi\in\H$, we have as $k\rightarrow\infty$,
$$Z_k(Hilb_k(\phi))\rightarrow E,$$
and the convergence is uniform in $C^l(l>>1)$ bounded subsets  of
$\H$.\\
\end{lem}
In section 3 we will investigate more about the relation
between $Z_k$ and $E$.\\

The following is proved essentially in $\cite{D5}$:\\

 \begin{lem}\label{convexity}(\cite{D5}) $\L_k$ is convex on $\H$, and
$Z_k$ is convex on $\B_k$.\\%
\end{lem}

\section{Convergence of geodesic distance}
In this section we shall prove theorem \ref{main}. The following
lemma about the derivative of $Hilb_k$ and $FS_k$ follows from a
simple
calculation:\\

\begin{lem} \label{tangentmap}For $\phi\in \H$, and $\delta\phi\in T_{\phi}\H$, we have for
any $s_1, s_2\in H^0(X, L^k)$, that
\begin{equation}
d _{\phi}\ Hilb_k(\delta\phi)(s_1, s_2)=\int (s_1,
s_2)_{\phi}(-k\delta\phi+\Delta\delta\phi)d\mu_{\phi}. \label{T1}
\end{equation}
For $H\in \B_k$, and $\delta H\in T_H \B_k$, we have
\begin{equation}
d_{H}\ FS_k(\delta H)=-\frac{1}{k}\sum_{i,j} \delta H(s_i,
s_j)\cdot(s_j, s_i)_{FS_k(H)}, \label{T2}
\end{equation}
where $\{s_i\}$ is an orthonormal basis of $H$.\\
\end{lem}

The following theorem proves the convergence of infinitesimal
geometry.\\

\begin{thm} \label{local convergence}Given a smooth path $\phi(t)\in \H(t\in [0,1])$, after
normalization, the length of the induced path $H_k(t)=Hilb_k
(\phi(t))$ converges to that of $\phi(t)$ as $k\rightarrow\infty$,
where the corresponding metrics on $\H$ and $\B_k$ are used. More
precisely,
$$\lim_{k\rightarrow\infty}k^{-n-2}||\dot{H}_k||^2=||\fid||^2=\int_X \dot{\phi}^2d\mu_{{\phi}},$$
for any $t\in[0,1].$\\
\end{thm}

\begin{proof} W.L.O.G, assume $t=0$.  Let $\{s_i(t)\}$ be an orthonormal
basis of $H_k(t)=Hilb_k(\phi(t))$. Then by lemma \ref{tangentmap},
$$\dot{H}_{ij}=\dot{H}(s_i(t), s_j(t))=\int_X ( s_i, s_j )(-k\dot{\phi}+\Delta\dot{\phi})d\mu_{\phi},$$
where for convenience we drop the subscript $k$ of  $H_k$.
Diagonalize it at $t=0$, we get
$$\dot{H}_{ij}=\delta_{ij}\cdot \int_X |s_i|^2_0(-k\dot{\phi}+\Delta\dot{\phi})d\mu_{\phi}.$$
Now let $\phi_k(t)=FS_k(H_k(t))$, then
$$\dot{\phi_k}=-\frac{1}{k^{}}\sum_i\dot{H}_{ii}|s_i|^2_{\phi_k}.$$
Therefore,
\begin{eqnarray*}
||\dot{H}||^2&=&\sum_i|\dot{H_{ii}}H_{ii}^{-1}|^2\\&=&  \int_X
\sum_i\dot{H}_{ii}|s_i|^2(-k\dot{\phi}+\Delta\dot{\phi})d\mu_{\phi}\\
  &=& -k^{} \int_X
\dot{\phi_k}(-k\dot{\phi}+\Delta\dot{\phi})\rho_k(\phi)d\mu_{\phi}\\
\end{eqnarray*}
By lemma \ref{tianyauzelditch}, we know
$$\rho_k(\phi)=k^n+\frac{1}{2}S\cdot k^{n-1}+O(k^{n-2}),$$ and
$$\dot{\phi_k}-\dot{\phi}=\frac{1}{k}\cdot\frac{\dot{\rho}_k(\phi)}{\rho_k(\phi)}=\frac{\dot{S}}{2k^2}+O(\frac{1}{k^3}).$$
Hence,
$$||\dot{H}||^2=k^{n+2}(\int_X\dot{\phi}^2d\mu_{\phi}+O(\frac{1}{k})). \ \
$$\\ \end{proof}

\begin{rmk} \label{important remark}Actually we have proved that for any smooth $\psi$, and $\{s_i\}$
an orthonormal basis of $H_k$,
$$\lim_{k\rightarrow\infty}k^{-n-2}\sum_{i, j} |\int_X(s_i, s_j )_{\phi}(-k\psi+\Delta\psi)d\mu_
{{\phi}}|^2=\int_X{\psi^2} d\mu_{{\phi}}.$$ Indeed, the
convergence is uniform for $\psi$ varying in a $C^l$ compact set.
So the following holds:
\begin{equation}
\lim_{k\rightarrow\infty}k^{-n}\sum_{i, j} |\int_X(s_i, s_j)
_{\phi}\psi d\mu_ {{\phi}}|^2=\int_X{\psi^2}
d\mu_{{\phi}}.\label{K}\end{equation}.\\
\end{rmk}

Now we can prove one side inequality of theorem \ref{main}.\\

\begin{cor}\label{first half}Given two metrics $\phi_1, \phi_2\in\H$, and denote
$$H_{k, i}=Hilb_k(\phi_i), \ \ i=1,2,$$then we have
$$\limsup_{k\rightarrow\infty}k^{-\frac{n}{2}-1}d_{\B_k}(H_{k,1}, H_{k,2})\leq
d_{\H}(\phi_1,\phi_2).$$
\end{cor}

\begin{proof} From lemma \ref{geodesic approximation lemma}, we know that
for any $\epsilon>0$, there exists a smooth $\epsilon$-approximate
geodesic $\phi(t)(t\in[0,1])$ in $\H$ connecting $\phi_1$ and
$\phi_2$ such that the length
$$L_{\H}(\phi(t))\leq d_{\H}(\phi_1, \phi_2)+\frac{\epsilon}{2}.$$
For this path and $k$ sufficiently large, by theorem \ref{local
convergence},
$$k^{-\frac{n}{2}-1}L_{\B_k}(Hilb_k(\phi(t)))\leq L_{\H}(\phi(t))+\frac{\epsilon}{2}.$$
Then $$k^{-\frac{n}{2}-1}d_{\B_k}(H_{k,1}, H_{k,2})\leq
k^{-\frac{n}{2}-1}L_{\B_k}(Hilb_k(\phi(t)))\leq
d_{\H}(\phi_1,\phi_2)+\epsilon. \ \ $$ \end{proof}

To prove the reversed inequality, we need the following lemma.\\

\begin{lem} \label{second derivative convergence} Given a smooth path $\phi(t)\in \H(t\in [0,1])$, then
$$\lim_{k\rightarrow\infty}k^{-n-2}||\nabla_{\dot{H}_k}\dot{H}_k||^2=
||\nabla_{\fid}\fid||^2=\int_X
|\nabla_{\dot{\phi}}\dot{\phi}|^2d\mu_{\phi},$$
for any $t\in[0,1].$\\
\end{lem}
\begin{proof} First notice that
$$\nabla_{\dot{\phi}}\dot{\phi}=\ddot{\phi}-|\nabla\dot{\phi}|^2,$$ and
$$\nabla_{\dot{H}}\dot{H}=\frac{d}{dt}(H^{-1}\dot{H})H=\ddot{H}-\dot{H}H^{-1}\dot{H},$$
where again we drop the subscript $k$ for convenience. As before,
we only need to prove the theorem at $t=0$. We can pick an
orthonormal basis $\{s_i(t)\}$ with respect to  $H(t)$ such that
$\dot{H}(0)$ is diagonalized, i.e.
\begin{equation*}
\int ( s_i(t), s_j(t))_{\phi} d\mu_{\phi}=\delta_{ij},
\end{equation*}
 and
\begin{equation}
\dot{H}_{ij}=\int (s_i, s_j)(-k\dot{\phi}+\Delta
\dot{\phi})d\mu=\dot{H}_{i}\cdot\delta_{ij},
\end{equation}
 holds at $t=0$. Taking more
derivatives, we obtain:
\begin{equation}
\ddot{H}_{ij}:=\ddot{H}(s_i, s_j)=\int (s_i, s_j)
[(-k\dot{\phi}+\Delta\dot{\phi})^2-k\ddot{\phi}
+\Delta{\ddot{\phi}}+\Psi(\dot{\phi})]d\mu,
\end{equation}
where $$\Psi(f)=-\sum_{i,j} f_{i\bar{j}}f_{j\bar{i}}.$$
\begin{eqnarray}
\dddot{H}_{ij}:=\dddot{H}(s_i, s_j)&=&\int ( s_i, s_j)
[(-k\fid+\Delta\fid)^3-k\fidd(-k\fid+\Delta\fid)+(\Delta\fidd+\Psi(\fid))(-k\fid+\Delta\fid)+2
k^2 \fid\fidd\nonumber\\&&-2k
\fidd\Delta{\fid}-2k\fid(\Delta{\fidd}+\Psi(\fid))+\frac{d}{dt}(\Delta{\fid})^2-k\fiddd+\frac{d^2}{dt^2}(\Delta\fid)]d\mu
\nonumber\\&=&\int(s_i, s_j)
[-k^3\fid^3+3k^2(\fid^2\Delta\fid+\fidd\fid)-3k(\fid(\Delta\fid)^2+\fidd\Delta{\fid}
+\fid\Delta\fidd+\fid\Psi(\fid)+\frac{1}{3}\fiddd)\nonumber\\&&+3\Delta\fid\frac{d}{dt}(\Delta\fid)+\frac{d^2}{dt^2}
(\Delta\fid)]d\mu.
\end{eqnarray}

Let $\dot{s}_i=\sum_{j}a_{ij}s_j$, then we can further assume that
$(a_{ij})$ is Hermitian, i.e $a_{ji}=\overline{a_{ij}}$. Then it
is easy to see that
$$a_{ij}=-\frac{1}{2}\dot{H}_{ij}.$$ Let $\phi_k(t)=FS_k(H_k(t))$, i.e.
$$\sum_{i}|s_i|^2_{h_0 e^{-k\phi_k}}=1,$$
and
\begin{equation}\dot{\phi_k}(t)=-\frac{1}{k^{}}\sum_{i,j}\dot{H}_{ij}
(s_j,
s_i)_{\phi_k}=-\frac{1}{k\rho_k(\phi)}\sum_{i,j}\dot{H}_{ij}( s_j,
s_i)_{\phi}. \end{equation} Taking time derivative, we get
\begin{equation}
\pidd(t)=\frac{\dot{\rho_k}}{k\rho_k^2}\sum_{i,
j}\dot{H}_{ij}(s_j,
s_i)_{\phi}-\frac{1}{k\rho_k}[\sum_{i,j}\ddot{H}_{ij}(s_j,
s_i)_{\phi}-2\sum_{i,j,l}\dot{H}_{il}\dot{H}_{lj}(s_j,
s_i)_{\phi}]-k\fid\pid,
\end{equation}
i.e.
\begin{equation}
\sum_{i,j}\ddot{H}_{ij}(s_j,
s_i)_{\phi}=-k\rho_k\pidd-k\dot{\rho_k}\pid -k^2 \rho_k\fid\pid+2
\sum_{i, j, l}\hd_{il}\hd_{lj}(s_j, s_i)_{\phi}.
\end{equation}
Define:
\begin{eqnarray*}
g_k&=&k^{-1}(-k+\Delta)^{-1}[(-k\dot{\phi}+\Delta\dot{\phi})^2-k\ddot{\phi}
+\Delta{\ddot\phi}+\Psi(\dot{\phi})]\\&=&-\fid^2-\frac{1}{k}(2|\nabla\fid|^2-\fidd)-
\frac{1}{k^2}[(\Delta\fid)^2+\Psi(\fid)-2\Delta(\fid\Delta\fid)+\Delta^2(\fid^2)]+O(\frac{1}{k^3}).
\end{eqnarray*}
Now let $\psi_k=d\ \text{FS}_k\circ d\ \text{Hilb}_k(g_k)$, then
\begin{eqnarray*}
-\frac{1}{k^2\rho_k}\sum_{i,j}\hdd_{ij}(s_j,
s_i)=\psi_k=g_k-\frac{\Phi(\fid^2)}{k^{2}}+O(k^{-3}),
\end{eqnarray*}
where $$\Phi(f)=\frac{1}{2}\delta_{f}S.$$ From (12) we then see
that
\begin{eqnarray*}
2 \sum_{i,j,l}\hd_{il}\hd_{lj}(s_j,
s_i)_{\phi}&=&-k^2\rho_kg_k+k\rho_k\pidd+k\dot{\rho_k}\pid+k^2\rho_k\fid\pid+\Phi(\fid^2)k^n+O(k^{n-1})\\
&=&2k^2\rho_k\fid^2+2k\rho_k|\nabla\fid|^2\\\\&&
+[(\Delta\fid)^2+\Psi(\fid)-2\Delta(\fid\Delta\fid)+2\fid\Phi(\fid)+\Phi(\fid^2)+\Delta^2(\fid^2)]
k^n+O(k^{n-1}).
\end{eqnarray*}
Thus,
\begin{eqnarray}
&&\sum_{i,j}(\hdd_{ij}-\hd_{i}\hd_{ij})(s_j, s_i)\nonumber\\
&=&k\rho_k(|\nabla\fid|^2-\fidd)+\frac{\rho_k}{2}[\Phi(\fid^2)+\Delta^2(\fid^2)
-2\fid\Phi(\fid)+(\Delta\fid)^2+\Psi(\fid)-2\Delta(\fid\Delta\fid)]+O(k^{n-1}).\nonumber\\
\end{eqnarray}
By (10), $$\hddd=k^2d\ \text{Hilb}_k(f_k),$$ where
\begin{eqnarray*}
f_k&=& k^{-2}(-k+\Delta)^{-1}[
-k^3\fid^3+3k^2(\fid^2\Delta\fid+\fidd\fid)-3k(\fid(\Delta\fid)^2+\fidd\Delta{\fid}
+\fid\Delta\fidd+\fid\Psi(\fid)+\frac{1}{3}\fiddd)\nonumber\\\\&&+3\Delta\fid\frac{d}{dt}(\Delta\fid)+\frac{d^2}{dt^2}
(\Delta\fid)]\nonumber\\\\
&=&
-\frac{1}{k^3}[-k^3\fid^3+3k^2(\fid^2\Delta\fid+\fidd\fid)-3k[\fid(\Delta\fid)^2+\fidd\Delta\fid+\fid\Delta\fidd
+\fid\Psi(\fid)+\frac{1}{3}\fiddd]\nonumber\\\\&&-(3k^2\fid^2\Delta\fid+6k^2\fid|\nabla\fid|^2)+3k[\Delta(\fid^2\Delta\fid)
+\Delta (\fidd\fid)-\Delta^2(\fid^3)]+O(1)] \nonumber\\\\&=&
\fid^3-\frac{3}{k}\fid(\fidd-2|\nabla\fid|^2)\\&&+\frac{3}{k^2}[\fid(\Delta\fid)^2+\fidd\Delta\fid+\fid\Delta\fidd+\fid\Psi(
\fid)+\frac{1}{3}\fiddd-\Delta(\fid^2\Delta\fid)-\Delta(\fidd\fid)-\Delta^2(\fid^3)]
+O(\frac{1}{k^3}).\nonumber\\\\\end{eqnarray*} Therefore,
\begin{eqnarray}
\sum_{i,j}\hddd_{ij}(s_j, s_i)&=&-k^3\rho_k
f_k-k^{n+1}\Phi(\fid^3)+O(k^n)\nonumber\\&=&
-k^3\rho_k\fid^3+3k^2\rho_k\fid(\fidd-2|\nabla\fid|^2)-k\rho_k[\Phi(\fid^3)
+3\fid(\Delta\fid)^2+3\fidd\Delta\fid+3\fid\Delta\fidd\nonumber\\\nonumber\\&&+3\fid\Psi(
\fid)+\fiddd-3\Delta(\fid^2\Delta\fid)-3\Delta(\fidd\fid)-\Delta^2(\fid^3)]+O(k^n).\nonumber\\
\end{eqnarray}
Now differentiating (13), we obtain
\begin{eqnarray}
&&\sum_{i,j}\hddd_{ij}(s_j,
s_i)-\sum_{i,j}\hdd_{i,j}(\hd_i+\hd_j)(s_j,
s_i)-k\fid\sum_{i,j}\hdd_{ij}(s_j,
s_i)\nonumber\\\nonumber\\&&=-2k
\dot{\rho_k}\pidd-k\rho_k\piddd-k\ddot{\rho_k}\pid-k^2\dot{\rho_k}\fid\pid-k^2\rho_k\fidd\pid-k^2\rho_k\fid\pidd
\nonumber\\\nonumber\\&&+2 \sum_{i,j}\hdd_{i,j}(\hd_i+\hd_j)(s_j,
s_i)-6\sum_{i} \hd_i^3|s_i|^2-2k\fid\sum_i \hd_i^2|s_i|^2.
\end{eqnarray}
Thus,
\begin{eqnarray}
&&3\sum_{i,j}\hdd_{i,j}(\hd_i+\hd_j)(s_j, s_i)-6\sum_i
\hd_i^3|s_i|^2\nonumber\\\nonumber\\&=&\sum_{i, j}\hddd_{i,
j}(s_j,
s_i)+k\fid(k\rho_k\pidd+k\dot{\rho_k}\pid+k^2\rho_k\fid\pid)\nonumber\\&&
+2k\dot{\rho_k}\pidd+k\rho_k\piddd+k\ddot{\rho_k}
\pid+k^2\dot{\rho_K}\fid\pid+k^2\rho_k\fidd\pid+k^2\rho_k\fid\pidd\nonumber\\\nonumber\\&=&\sum_{i,
j}\hddd_{i,j}(s_j,
s_i)+k^3\rho_k\fid^2\pid+3k^2\rho_k\fid\fidd+k^{n+1}(2\Phi(\fid)\fid^2+\fiddd)+O(k^n)\nonumber\\\nonumber\\&=&
6k^2\rho_k\fid(\fidd-|\nabla\fid|^2)+k^{n+1}[3\Phi(\fid)\fid^2+\fiddd-\Phi(\fid^3)-3\fid(\Delta\fid)^2-3\fid\Delta\fidd
-3\fidd\Delta\fid\nonumber\\\nonumber\\&&-3\fid\Psi(\fid)-\fiddd+3\Delta(\fid^2\Delta\fid)+3\Delta(\fidd\fid)-\Delta^2
(\fid^3)]+O(k^n) \nonumber\\\nonumber\\&=&
6k^2\rho_k\fid(\fidd-|\nabla\fid|^2)+k^{n+1}[3\Phi(\fid)\fid^2-\Phi(\fid^3)-\Delta^2(\fid^3)-3\fid(\Delta\fid)^2-6\nabla\fid
\cdot\nabla\fidd
\nonumber\\\nonumber\\&&-3\fid\Psi(\fid)+3\Delta(\fid^2\Delta\fid)]+O(k^n).
\end{eqnarray}
Putting (8), (9), (14), (17) together, we get
\begin{eqnarray*}
|\nabla_{\hd}\hd|^2&=&\sum_{i,j}(\hdd_{ij}-\hd_i\hd_{ij})\int(s_j,
s_i)(k^2\fid^2-2k\fid\Delta\fid-k\fidd+O(1))d\mu\\\\&-&\sum_{i,j}(\hdd_{ij}\hd_i-\hd_i^3)\int(s_j,
s_i)(-k\fid+\Delta\fid)d\mu\\\\&=&
\int\{k\rho_k(|\nabla\fid|^2-\fidd)+\frac{\rho_k}{2}[\Phi(\fid^2)
-2\fid\Phi(\fid)+\Delta^2(\fid^2)+(\Delta\fid)^2+\Psi(\fid)-2\Delta(\fid\Delta\fid)]+O(k^{n-1})\}\\\\&&
(k^2\fid^2-2k\fid\Delta\fid-k\fidd+O(1))d\mu\\\\&&+\int
\{k^2\rho_k\fid(\fidd-|\nabla\fid|^2)+\frac{k\rho_k}{6}[3\Phi(\fid)\fid^2-\Phi(\fid^3)-\Delta^2(\fid^3)\\\\&&+
6\nabla\fid\nabla\fidd-3\fid\Psi(\fid)-3\fid(\Delta\fid)^2+3\Delta(\fid^2\Delta\fid)+O(k^n)]\}(k\fid-\Delta\fid)d\mu
\\\\&=&\int
-k^2\rho_k(|\nabla\fid|^2-\fidd)(2\fid\Delta\fid+\fidd)\\\\&&
+\frac{k^2}{2}\rho_k\fid^2[\Phi(\fid^2)
-2\fid\Phi(\fid)+\Delta^2(\fid^2)+(\Delta\fid)^2+\Psi(\fid)-2\Delta(\fid\Delta\fid)]\\\\&&-k^2\rho_k\fid\Delta\fid(\fidd-|\nabla\fid|^2)
\\\\&&
+\frac{k^2\rho_k}{6}\fid[3\Phi(\fid)\fid^2-\Phi(\fid^3)-\Delta^2(\fid^3)+
6\nabla\fid\nabla\fidd-3\fid\Psi(\fid)-3\fid(\Delta\fid)^2+3\Delta(\fid^2\Delta\fid)]d\mu+O(k^{n+1})
\\\\&=&-\frac{k^2}{2}\int \rho_k\Delta(\fid^2)|\nabla\fid|^2d\mu
+\frac{k^2}{2}\int\rho_k\fid^2[-R_{i\bar{j}}\fid_{\bar{i}}\fid_j+\Delta^2(\fid^2)+\fid\Delta^2\fid
+(\Delta\fid)^2
-2\Delta(\fid\Delta\fid)]d\mu\\\\&&+\frac{k^{n+2}}{6}\int\fid[3R_{i\bar{j}}\fid\fid_{\bar{i}}\fid_j-\frac{3}{2}\fid^2\Delta^2
\fid+\frac{1}{2}\Delta^2(\fid^3)-\Delta^2(\fid^3)-3\fid(\Delta\fid)^2+3\Delta(\fid^2\Delta\fid)]d\mu+O(k^{n+1}).\\\\&=&
-k^2\int
\rho_k(\fid\Delta\fid+|\nabla\fid|^2)|\nabla\fid|^2d\mu\\\\&&+k^{n+2}\int
\frac{1}{2}\fid^3\Delta^2\fid+\frac{1}{4}\fid^2\Delta^2(\fid^2)+\frac{1}{2}\fid^2(\Delta\fid)^2-2\fid^2(\Delta\fid)^2-2\fid\Delta\fid|\nabla\fid|^2
d\mu\\\\&& -\frac{k^{n+2}}{3}\int \fid^3\Delta^2\fid
d\mu+O(k^{n+1})\\\\&=&k^{n+2}\int(\fidd-|\nabla\fid|)^2d\mu+O(k^{n+1}).\\\\
\end{eqnarray*}
where we used:
$$\Phi(\fid)=-\frac{1}{2}(\fid_{i\bar{j}}R_{j\bar{i}}+\Delta^2\fid),$$
$$\Phi(\fid^2)=-\fid\fid_{i\bar{j}}R_{j\bar{i}}-\fid_i\fid_{\bar{j}}R_{j\bar{i}}-\frac{1}{2}\Delta^2(\fid^2),$$
$$\Phi(\fid^3)=-\frac{3}{2}\fid^2\fid_{i\bar{j}}R_{j\bar{i}}-3\fid\fid_{i}\fid_{\bar{j}}R_{j\bar{i}}-
\frac{1}{2}\Delta^2(\fid^3).$$
$$\int \fid^2\Delta^2(\fid^2)d\mu=4\int\fid^2(\Delta\fid^2)d\mu+8\int \fid\Delta\fid|\nabla\fid|^2d\mu+4\int |\nabla\fid|^4d\mu.$$
$$\int \fid^3\Delta^2(\fid)d\mu=3\int \fid^2(\Delta\fid)^2d\mu+6\int
\fid\Delta\fid|\nabla\fid|^2d\mu.\ \ \ \ $$
\\\\
\end{proof}
 \begin{proof} \emph{of Theorem \ref{main}}. In view of corollary \ref{first half}, it
suffices to show that
$$\liminf_{k\rightarrow\infty}k^{-\frac{n}{2}-1}d_{\B_k}(H_{k,0}, H_{k,1})\geq
d_{\H}(\phi_0,\phi_1).$$ By lemma \ref{geodesic approximation
lemma}, for any $\epsilon>0$ small enough, there exists a smooth
family $\phi_{\epsilon}(\cdot):[0,1]\rightarrow\H$, such that
$\phi_{\epsilon}(0)=\phi_0$, $\phi_{\epsilon}(1)=\phi_1$, and
$$(\fidd_{\epsilon}-|\nabla\fid_{\epsilon}|_t^2)\omega_{\epsilon,t}^n=\epsilon\omega^n.$$
Moreover, $\fidd_{\epsilon}-|\nabla\fid_{\epsilon}|_t^2\leq C$,
where $C$ is a constant independent of $t$ and $\epsilon$. Let
$H_{\epsilon}(t)=Hilb_k(\phi_{\epsilon}(t))$, then by lemma
\ref{second derivative convergence}, for $k$ large enough, we have
$$k^{-n-2}|\nabla_{\hd_{\epsilon}}\hd_{\epsilon}|^2\leq C\epsilon.$$ By the following
simple lemma, we know that
$$k^{-\frac{n}{2}-1}d_{\B_k}(H_{k,0}, H_{k,1})\geq k^{-\frac{n}{2}-1}
L(H_{\epsilon})-\sqrt{C\epsilon}\rightarrow
L(\phi_{\epsilon})-\sqrt{C\epsilon},$$ as $k\rightarrow \infty$.
Hence,
$$\liminf_{k\rightarrow\infty}k^{-\frac{n}{2}-1}d_{\B_k}(H_{k,0},
H_{k,1})\geq d_{\H}(\phi_0,\phi_1)-\sqrt{C\epsilon}.$$ Let
$\epsilon\rightarrow 0$, we obtain the desired result. \\
\end{proof}
Now we prove a simple lemma from Riemannian geometry.\\

\begin{lem}\label{Riemannian lemma 1}Suppose $(M, g)$ is a simply connected Riemmannian manifold
with non-positive curvature. Let $\gamma:[0,1]\rightarrow M$ be a
path in $M$ such that
$|\nabla_{\dot{\gamma}(t)}\dot{\gamma}(t)|\leq \epsilon$ for all
$t\in[0,1]$. Then
$$d(\gamma(0), \gamma(1))\geq L(\gamma)-\epsilon.$$
\end{lem}

\begin{proof} Denote $p=\gamma(0)$ and $q=\gamma(1)$. By the theorem of
Cartan-Hadamard, we know that the exponential map at $p$ is a
diffeomorphism. Now let $\gamma_t(s)$$(s\in [0,1])$ be the unique
geodesic connecting $p$ and $\gamma(t)$. By standard calculation
of the second variation, we obtain:
\begin{eqnarray*}
\frac{d^2}{dt^2}L(\gamma_t)&=&\frac{1}{d_t}\langle\nabla_{\dot{\gamma}(t)}\dot{\gamma}(t),
\frac{\partial} {\partial s}\gamma_t(s)|_{s=1}\rangle
\\&&+\frac{1}{d_t}\int_0^1|(\frac{\p^2}{\p s\p
t}\gamma_t(s))^{\perp}|^2-R(\frac{\p}{\p s}\gamma_t(s),
\frac{\p}{\p t}\gamma_t(s), \frac{\p}{\p s}\gamma_t(s),
\frac{\p}{\p t}\gamma_t(s) )ds \\&& \geq -\epsilon,
\end{eqnarray*}
where $d_t=L(\gamma_t)$, and $\perp$ denotes projection to the
orthogonal of $\frac{\p}{\p s}\gamma_t(s)$. It is also easy to see
that
$$L(\gamma_0)=0,$$
and $$\frac{d}{dt} L(t)|_{t=0}=|\dot{\gamma}(0)|.$$ Therefore,
$$L(t)\geq |\dot{\gamma}(0)|t-\frac{\epsilon}{2}t^2.$$
In particular, $d(p, q)=L(\gamma(1))\geq
 |\dot{\gamma}(0)|-\frac{\epsilon}{2}$. On the other hand,
$$L(\gamma)=\int_0^1|\dot{\gamma}(t)|dt\leq \int_0^1 (|\dot{\gamma}(0)|+\epsilon t) dt=|\dot{\gamma}(0)|+\frac{\epsilon}{2}.$$
Hence, $d(p,q)\geq L(\gamma)-\epsilon$.
\end{proof}

\section{Weak convexity of K-energy}
In this section, we shall prove corollary \ref{inequality 1} and
\ref{inequality 2}. Before doing this, we want to show the
functional $I$ is quantized by
$I_k$, which is an analogue of the fact that $Z_k$ quantizes $E$.\\

\begin{prop} There are constants $c_k$ such that for $\phi\in \H$,
and $H_k=Hilb_k(\phi)$, we have
$$\lim_{k\rightarrow\infty}k^{-n-1}I_k(H_k)+c_k=I(\phi).$$ The convergence is uniform when $\phi$ varies in a
 $C^l$ bounded sets as before.\\
\end{prop}

\begin{proof}It suffices to show for any $\psi\in T_{\phi}\H$
\begin{equation}
d_{\phi}\ I(\psi)=d_{Hilb_k(\phi)}\ I_k\circ d_{\phi}\
Hilb_k(\psi). \label{I}
\end{equation}
By definition,
$$d_{\phi}\ I(\psi)=\int \psi d\mu_{\phi}.$$
On the other hand, by lemma \ref{tangentmap},
\begin{eqnarray*}
d_{Hilb_k(\phi)}\ I_k\circ d_{\phi}\ Hilb_k(\psi)&=&\sum_i\int
|s_i|_{\phi}^2 (-k\psi+\Delta\psi)d\mu_{\phi}\\ &=& \int
\rho_k(\phi)(-k\psi+\Delta\psi)d\mu_{\phi}
\end{eqnarray*}
Then the result follows from lemma \ref{tangentmap}.
\end{proof}

\begin{prop} \label{gradient convergence} Let $\phi\in \H$, and $H_k=Hilb_k(\phi)$, then
$$\lim_{k\rightarrow\infty}k^{n+2}||
\nabla Z_k(H_k)||^2=||\nabla E(\phi)||^2.$$
\end{prop}

\begin{proof} Note $$||\nabla E(\phi)||^2=\int
(S-\underline{S})^2d\mu_{\phi}$$ is simply the Calabi's
functional. An inequality of this form is essentially proved in
\cite{D6} for obtaining a lower bound of the Calabi functional. We
can easily calculate the first variation of $Z_k$:
$$\delta Z_k=-\frac{ d_k}{k^n\cdot V}\int_X \delta H_{ij}[( s_j, s_i )_{FS_k(H)}]_0
d\mu_{\omega_{FS_k(H)}}, $$ where $[A]_0$ denote the trace-free
part of a matrix $A$, and $\{s_i\}$ is an orthonormal basis with
respect to $H$. Here
$$V=\int_X d\mu_{\phi}, $$ and $$d_k=dim H^0(X, L^k)=\int_X
\rho_k(\phi)d\mu_{\phi}=k^n\cdot V+\frac{k^{n-1}\cdot
V}{2}\underline{S}+O(k^{n-2}).$$ So
$$(\nabla Z_k)_{ij}=-\frac{  d_k}{k^n\cdot V} \int_X [( s_i, s_j )_{FS_k(H)}]_0 d\mu_{\omega_{FS_k(H)}}.$$
Diagonalize the above matrix so that its diagonal entries are
$$\lambda_{i}=\frac{ d_k}{k^n\cdot V} \int_X
(|s_i|^2_{\psi_2}-\frac{1}{d_k}) d\mu_{{\psi_2}}.$$ Let
$\phi_k=FS_k(H_k)$, then
$$\phi_k-\phi=\frac{1}{k}\log(\rho_k(\phi))=\frac{1}{k}\log(k^n+\frac{S}{2 }k^{n-1}+O(k^{n-2})).$$
So
$$|s_i|_{\phi_k}^2=k^{-n}|s_i|_{\phi}^2(1-\frac{S}{2 k}+O(\frac{1}{k^2})),$$
and
$$\omega_{\phi_k}^n=\omega_{\phi}^n(1+O(k^{-2})).$$
Since $$\int_X |s_i|_{\phi}^2d\mu_{{\phi}}=1,$$ we obtain
$$\lambda_i=-\frac{1}{k^{n+1}}\int_X |s_i|_{\phi}^2(S-\underline{S})d\mu_{{\phi}}+O(\frac{1}{k^{n+2}}).$$
Hence by remark \ref{important remark},
\begin{eqnarray*}
||\nabla Z_k||^2&=&\sum_i|\lambda_i|^2\\
&\leq& k^{-2n-2}\sum_i[\int_X
|s_i|_{\phi}^2(S-\underline{S})d\mu_{{\phi}}+O(\frac{1}{k^{n+2}})]^2\\&=&
k^{-n-2}(Ca(\phi_2)+O(\frac{1}{k})). \ \ \ \
\end{eqnarray*}

\end{proof}
The proof of the following proposition is similar and we omit
it.\\

 \begin{prop}\label{derivative convergence}Let $\phi(t)$ be a smooth path
in $\H$, and $H_k(t)=Hilb_k(\phi)$, then
$$\lim_{k\rightarrow\infty}\frac{d}{dt}Z_k(H_k)=\frac{d}{dt}E(\phi).$$\\
\end{prop}

 \begin{proof}  \emph{of  corollary \ref{inequality 1}}. By lemma \ref{convexity}, we already know that $Z_k$ is a
genuine convex function on $\B_k$. So it is easy to see that
$$Z_k(H_2)-Z_k(H_1)\leq d_{\B_k}(H_1, H_2)\cdot ||\nabla Z_k(H_2)||.$$
Using theorem \ref{main}, lemma \ref{ZquantizesE} and proposition
\ref{gradient convergence}, and let $k\rightarrow\infty$, we get
the desired inequality for $E$. \ \ \\
\end{proof}

\begin{proof}\emph{of corollary \ref{inequality 2}.} As in the proof of
theorem \ref{main}, for any $\epsilon>0$ sufficiently small, we
choose $\epsilon$-approximate geodesic
$\phi_{\epsilon}:[0,1]\rightarrow \H$, such that
$|\nabla_{\fid_{\epsilon}}\fid_{\epsilon}|\leq C\epsilon$. Denote
$H_{\epsilon}(t)=Hilb_k(\phi_{\epsilon}(t))$,
%$\psi_{\epsilon}(t)=FS_k(H_{\epsilon}(t))$,
and $H_i=Hilb_k(\phi_i)$, $i=0,1$. Then by lemma \ref{second
derivative convergence}, we see that
$$k^{-\frac{n}{2}-1}|\nabla_{\hd_{\epsilon}}\hd_{\epsilon}|\leq C\epsilon,$$
for $k$ sufficiently large and a different constant $C$. Denote by
$\tilde{H}:[0, 1]\rightarrow\B_k$ the geodesic connecting $H_0$
and $H_1$. By lemma \ref{Riemannian lemma 2} which will be proved
later, we know that
$$k^{-\frac{n}{2}-1}|\hd_{\epsilon}(i)-\dot{\tilde{H}}(i)|\leq
C\sqrt{\epsilon},$$ for $i=0,1$. Since $\tilde{H}(t)$ is a
geodesic in $\B_k$, we have by lemma \ref{convexity} that
$$dZ_{k_{H_1}}(\dot{\tilde{H}}(1))\geq dZ_{k_{H_0}}(\dot{\tilde{H}}(0)).$$
On the other hand, by proposition \ref{gradient convergence},
$$|dZ_{k_{H_1}}(\dot{H}_{\epsilon}(1)-\dot{\tilde{H}}(1))|\leq |\nabla Z_k|_{H_1}\cdot
|\dot{H}_{\epsilon}(1)-\dot{\tilde{H}}(1)|\leq
C\sqrt{\epsilon}(\sqrt{Ca(\phi_1)}+O(\frac{1}{k})).$$ Similarly,
$$|dZ_{k_{H_0}}(\dot{H}_{\epsilon}(0)-\dot{\tilde{H}}(0))|\leq C\sqrt{\epsilon}(\sqrt{Ca(\phi_0)}+O(\frac{1}{k})).$$
So,$$dZ_{k_{H_1}}(\dot{{H}}_{\epsilon}(1))\geq
dZ_{k_{H_0}}(\dot{{H}}_{\epsilon}(0))-C\sqrt{\epsilon}.$$ By
proposition \ref{derivative convergence}, that as $k\rightarrow
\infty$
$$dZ_{k_{H_1}}(\dot{{H}}_{\epsilon}(1))\rightarrow
dE_{\phi_1}(\dot{\phi}_{\epsilon}(1)),$$and
$$dZ_{k_{H_0}}(\dot{{H}}_{\epsilon}(0))\rightarrow
dE_{\phi_0}(\dot{\phi}_{\epsilon}(0)).$$ Thus,
$$dE_{\phi_1}(\dot{\phi}_{\epsilon}(1))\geq dE_{\phi_0}(\dot{\phi}_{\epsilon}(0))-C\sqrt{\epsilon}.$$
Letting $\epsilon\rightarrow 0$, this proves corollary
\ref{inequality 2}.
\\
\end{proof}

\begin{rmk} If we denote $\phi_{k,\epsilon}(t)=FS_k(Hilb_k{\phi_{\epsilon}(t)})$
 and $\tilde{\phi}_k(t)=FS_k(H(t))$, then by a similar estimate as
 in the proof of proposition \ref{gradient convergence}
 we have $$|\dot{\tilde{\phi}}_k(0)-\dot{\phi}_{k,\epsilon}(0)|_{L^2}\leq
 C\sqrt{\epsilon}.$$ Thus, the time derivative of $\tilde{\phi}_k$ at the
 end points converge in $L^2$ to the time derivative of the
 $C^{1,1}$ geodesic. This seems to be an interesting fact, compare
 \cite{PS1}.\\
\end{rmk}

 \begin{lem} \label{Riemannian lemma 2}Suppose $(M, g)$ is a simply connected Riemmannian
manifold with non-positive curvature. Let $\gamma:[0,1]\rightarrow
M$ be a path in $M$ such that
$|\nabla_{\dot{\gamma}(t)}\dot{\gamma}(t)|\leq \epsilon$ for all
$t\in[0,1]$. Let $\tilde{\gamma}:[0,1 ]$ be the unique geodesic
segment joining $\gamma(0)$ and $\gamma(1)$. Then
$$|\dot{\gamma}(0)-\dot{\tilde{\gamma}}(0)|^2\leq \frac{9}{4}\epsilon^2+4\epsilon
|\dot{\gamma}(1)|,$$ and
$$|\dot{\gamma}(1)-\dot{\tilde{\gamma}}(1)|^2\leq \frac{9}{4}\epsilon^2+4\epsilon
|\dot{\gamma}(0)|,$$
\end{lem}

\begin{proof} It suffices to prove the second inequality. As before, let
$\gamma_t(s):[0, 1]\rightarrow M$ be the smooth family of geodesic
segments connecting $\gamma(0)$ and $\gamma(t)$. Then by the proof
of lemma \ref{Riemannian lemma 1} we see that
$$\frac{d}{d t}\langle\frac{d\gamma}{d t},
\frac{\dot{\gamma}_t(1)}{|\dot{\gamma}_t(1)|}\rangle\geq
-\epsilon, $$ and so
$$\langle\frac{d\gamma}{d t},
\frac{\dot{\gamma}_t(1)}{|\dot{\gamma}_t(1)|}\rangle\geq
\lim_{t\rightarrow 0}\langle\frac{d\gamma}{d t},
\frac{\dot{\gamma}_t(1)}{|\dot{\gamma}_t(1)|}\rangle-\epsilon
t=|\frac{d\gamma}{d t}(0)|-\epsilon t.$$ Let
$$\frac{d\gamma}{d t}=A_t\cdot \frac{\dot{\gamma}_t(1)}{|\dot{\gamma}_t(1)|}+B_t$$
be the orthogonal decomposition. Then
$$A_t\geq |\frac{d\gamma}{d t}(0)|-\epsilon t.$$
It is clear that $$|\frac{d}{d t }|\frac{d\gamma}{d t}||\leq
|\frac{d ^2\gamma}{d t^2}|\leq \epsilon,$$ so
$$|\frac{d\gamma}{d t}(0)|-\epsilon t\leq|\frac{d\gamma}{d t}|\leq |\frac{d\gamma}{d t}(0)|+\epsilon t.$$
Hence,
$$|B_t|^2\leq(|\frac{d\gamma}{d t}(0)|+\epsilon t)^2-(|\frac{d\gamma}{d t}(0)|-\epsilon t)^2=4\epsilon t |\frac{d\gamma
}{dt}(0)|. $$ From the proof of lemma \ref{Riemannian lemma 1}, we
know that
$$|\frac{d\gamma}{dt}(0)|t-\frac{\epsilon}{2}t^2\leq|
\dot{\gamma}_t(1)|\leq
|\frac{d\gamma}{dt}(0)|t+\frac{\epsilon}{2}t^2.$$Finally,
\begin{eqnarray*}
|\frac{d\gamma}{dt}-\frac{1}{t}\dot{\gamma}_t(1)|^2&=&(A_t-\frac{1}{t}|\dot{\gamma}_t(1)|)^2+B_t^2\\
&\leq&\frac{9}{4}\epsilon^2t^2+4\epsilon t|\frac{d\gamma}{dt}(0)|,
\end{eqnarray*}
Let $t=1$, we obtain the desired bound.
\end{proof}

\section{Open problems}
In this section, we list some open problems and speculations on
which the results of this paper may help in the future.

\begin{q}
Theorem \ref{main} essentially says that $\H$ is a weak
Gromov-Hausdorff limit of a sequence of finite dimensional
symmetric spaces and $\H$ can be viewed as a generalized $Cat(0)$
space which has been extensively studied in the literature(c.f.
\cite{BH}). It would be interesting to investigate this more
carefully and develop a suitable notion for this type of
convergence. From this convergence it follows that the negativity
of curvature is inherited by the limit. A natural question would
be what else properties for the finite dimensional symmetric
spaces would survive on $\H$. For example can we describe the
algebraic structure of $\H$. In particular does it admit a local
involution? In other words, for any $\phi \in \cal H$, there
exists a small constant $\delta(\phi)$  such that there is an
local involution
\[
\sigma: B_\delta(\varphi) \bigcap {\cal H} \rightarrow
B_\delta(\varphi) \bigcap {\cal H}
\]
with $\sigma^2 = id.$ This is not trivial since the initial value
problem for the geodesic equation in $\H$ is not well-posed.

\end{q}

We know that $\cal H$ is not complete, then the following is very
interesting:

\begin{q} What is the structure of $\partial {\cal H}$?
\end{q}
It would be extremely interesting to understand the structure of
boundary. Note that for every $k$, the Bergman metric space $\B_k$
is complete and $\H$ is a limit of these nice symmetric spaces
after appropriate scaling down(by a factor of
$k^{-\frac{n+2}{2}}$). It is natural to hope that $\H$ should look
like the limit of the tangent cone of $\B_k$ at infinity.\\

By the work of J. Fine \cite{Fi}, we can approximate Calabi flow
over a bounded interval by the Balancing flow in $\B_k$. The next
intriguing question is can we approximate the Calabi flow over
$[0,\infty)$ if we know it converges (the  complex structure might
jump). This is related to

\begin{q} In an algebraic manifold, if there is a cscK metric, does the Calabi flow
exist for long time and converge to a cscK metric?
\end{q}

This is only proved for metrics near a cscK metric.  The
corresponding problems for
K\"ahler-Ricci flow was proved first by G. Perelman and a written version is provided by Tian-Zhu(\cite{TZ}).\\

\begin{q} The existence of cscK metrics implies that the K energy is
proper.
\end{q}

The corresponding results for K\"ahler-Einstein metrics is proved
by Tian in \cite{T}. Note the definition of properness could vary.
In \cite{T}, properness means bounded below by another positive
functional $J$. Since the K energy is convex along smooth
geodesics and its Hessian at a cscK metric is strictly positive
(if the automorphism group of the manifold is discrete), we would
expect it bounds the distance function. Note the corresponding
statements in the finite dimensional case is clear: a convex
function with a non-degenerate critical point automatically bounds
the distance function. Since both the distance function and the K
energy could be approximated by the corresponding  quantities on
$\B_k$, to prove the corresponding statement for $\H$, it is
important to derive
a uniform constant.\\

Xiuxiong Chen, Department of Mathematics, University of
Wisconsin-Madison, 480 Lincoln Drive, Madison, WI 53706.\\ Email:
xxchen@math.wisc.edu\\\\
 Song Sun, Department of Mathematics, University of Wisconsin-Madison,
 480 Lincoln Drive, Madison, WI 53706.\\ Email:
ssun@math.wisc.edu\\\\

\vskip3mm

\end{document}